\documentclass[11pt]{amsart}

\usepackage{amsmath}
\usepackage{mathrsfs}
\usepackage{amssymb}
\usepackage{tikz}
\usetikzlibrary{matrix,backgrounds}
\usepackage{comment}
\usepackage{color}
\usepackage[colorlinks,citecolor=blue,urlcolor=black,linkcolor=black]{hyperref}

\newtheorem{theorem}{Theorem}[section]
\newtheorem{claim}[theorem]{Claim}

\newtheorem{lemma}[theorem]{Lemma}

\theoremstyle{definition}
\newtheorem{definition}[theorem]{Definition}

\newtheorem{question}[theorem]{Question}

\theoremstyle{remark}

\newcount\skewfactor
\def\mathunderaccent#1#2 {\let\theaccent#1\skewfactor#2
\mathpalette\putaccentunder}
\def\putaccentunder#1#2{\oalign{$#1#2$\crcr\hidewidth
\vbox to.2ex{\hbox{$#1\skew\skewfactor\theaccent{}$}\vss}\hidewidth}}

\def\smallbox#1{\leavevmode\thinspace\hbox{\vrule\vtop{\vbox
   {\hrule\kern1pt\hbox{\vphantom{\tt/}\thinspace{\tt#1}\thinspace}}
   \kern1pt\hrule}\vrule}\thinspace}

\newcommand{\cf}{{\rm cf}}

\def\qedref#1{$\qed_{\reforiginal{#1}}$}

\title{An almost strong relation}
\author{Shimon Garti}
\address{Institute of Mathematics,
 The Hebrew University of Jerusalem,
 Jerusalem 91904, Israel}
\email{shimon.garty@mail.huji.ac.il}

\author{Andr\'es Villaveces}
\address{Departamento de Matem\'aticas, Universidad Nacional de Colombia, Bogot\'a 111321, Colombia}
\email{avillavecesn@unal.edu.co}


\begin{document}
\let\labeloriginal\label
\let\reforiginal\ref

\begin{abstract}
Let $\mu$ be a strong limit singular cardinal.
We prove that if $2^\mu>\mu^+$ then $\binom{\mu^+}{\mu}\rightarrow \binom{\tau}{\mu}_{<\cf(\mu)}$ for every ordinal $\tau<\mu^+$.
We obtain an optimal positive relation under $2^\mu=\mu^+$, as after collapsing $2^\mu$ to $\mu^+$ this positive relation is preserved.
\end{abstract}

\maketitle

\newpage

\section{Introduction}

The \emph{strong polarized partition relation}
$\binom{\lambda}{\kappa}\rightarrow \binom{\lambda}{\kappa}_\theta$ says that for every coloring
$c:\lambda\times\kappa\rightarrow\theta$ one can find $A\in[\lambda]^\lambda, B\in[\kappa]^\kappa$ such that $c''(A\times{B})$ is constant.
The \emph{almost strong polarized partition relation}
$\binom{\lambda}{\kappa}\rightarrow \binom{\tau}{\kappa}_\theta$ for
every $\tau<\lambda$ asserts that for every coloring
$c:\lambda\times\kappa\rightarrow\theta$ and any ordinal
$\tau<\lambda$ one can find a color $i<\theta$ and a pair of sets
$A\subseteq\lambda, B\subseteq\kappa$ such that ${\rm otp}(A)=\tau,
|B|=\kappa$ and $c''(A\times B)=\{i\}$.
The purpose of this paper is to prove such a relation at strong limit
singular cardinals.
The main theorem of the paper improves both a result of
Erd\H{o}s-Hajnal-Rado for strong limit singular cardinals with
countable cofinality and a result of Shelah for strong limit singular
cardinals with uncountable cofinality.

Erd\H{o}s, Hajnal and Rado showed in Theorem 42 of \cite{MR0202613} that whenever $\cf(\mu) =
\omega$, $\binom{\mu^+}{\mu} \rightarrow \binom{\mu}{\mu}_2$.
Assuming that $2^\mu>\mu^+$ and $\mu$ is a strong limit cardinal with countable cofinality, we
increase the order type of the first coordinate, namely ${\rm otp}(A)=\tau$ for every ordinal
$\tau<\mu^+$. 
Similarly, Shelah proved in \cite{MR1606515} that $\binom{\mu^+}{\mu} \rightarrow \binom{\mu+1}{\mu}_{<\cf(\mu)}$ whenever $\mu>\cf(\mu)>\omega$, $\mu$ is strong limit and $2^\mu>\mu^+$.
Actually, Shelah's theorem applies to singular cardinals of countable cofinality as well, but in
these cardinals the partition relation is less interesting since it holds even if without the assumption that $2^\mu>\mu^+$.
Anyway, we show that the ordinal $\mu+1$ can be replaced by every
ordinal $\tau<\mu^+$ for singular cardinals of uncountable cofinality under the same assumptions of $2^\mu>\mu^+$.
These results give a positive answer to Question 4.8 from \cite{MR3509813}.
Our proof follows in the footsteps of Shelah's proof in
\cite{MR1606515}, but we add a new feature which enables
us to stretch monochromatic sets of size $\mu$ and to obtain sets of
the same cardinality but a larger order type.

It is worth noting that the strong relation $\binom{\mu^+}{\mu} \rightarrow
\binom{\mu^+}{\mu}_2$ is consistent for singular cardinals (see
\cite{MR2987137}, \cite{MR3509813}), using several assumptions about
the structure of cardinal arithmetic.
For this relation the assumption $2^\mu>\mu^+$ is necessary since $2^\mu=\mu^+$ implies $\binom{\mu^+}{\mu} \nrightarrow
\binom{\mu^+}{\mu}_2$ as proved in \cite{MR0202613}, and usually the positive strong relation requires more than just $2^\mu>\mu^+$.
But for our result, it suffices to assume $2^\mu>\mu^+$ in
order to obtain the almost strong relation at strong limit singular
cardinals.
Moreover, although collapsing $2^\mu$ to $\mu^+$ destroys the strong
relation $\binom{\mu^+}{\mu} \rightarrow \binom{\mu^+}{\mu}_2$ it
still preserves the almost strong relation, as will be indicated at
the end of the paper.

The result of Shelah was introduced in an expository article of Menachem
Kojman, \cite{MR1356541}, who simplified the proof. The presentation here owes a lot to Kojman's paper, and in particular we follow the
notation of that paper.
We divide the rest of the paper into two sections.
In the first one we discuss the general concept of pcf arrays of
elementary submodels.
In the second section we prove the combinatorial theorem.
The main reason for this separation is that pcf arrays seem to
represent a general method.
Shelah remarked that tentatively this concept might lead to other
mathematical proofs, though no such one has been discovered so far.

Our notation is standard.
When we use elementary submodels of some $\mathcal{H}(\chi)$ we mean
that $\chi$ is a sufficiently large regular cardinal and this
structure can be augmented by any finite number of additional
predicates (like a well-ordering of the elements of
$\mathcal{H}(\chi)$).
If $\kappa=\cf(\kappa)<\lambda$ then $S^\lambda_\kappa =
\{\delta\in\lambda: \cf(\delta)=\kappa\}$. We shall use this notation
in most cases when $\lambda$ is a regular cardinal, in which case
$S^\lambda_\kappa$ is a stationary subset of $\lambda$.
An ordinal $\eta$ is indecomposable iff it has the form
$\omega^\beta$, where this denotes ordinal
exponentiation. Indecomposable ordinals behave like cardinals in the
sense that if $\alpha<\eta$ then the order type of $\eta - \alpha$ is
still $\eta$. The collection of all indecomposable ordinals of some
$\lambda=\cf(\lambda)>\aleph_0$ is a club subset of $\lambda$.
If $E$ is a club subset of $\lambda$ then ${\rm acc}(E)$ is the set of
accumulation points of $E$, and it is a club of $\lambda$ as well.

For basic background in pcf theory we suggest \cite{MR2768693} and \cite{MR1086455}.
For advanced theorems, including the main tool used in this paper, we
refer to \cite{MR1318912}.
A good background about polarized relations
appears in \cite{MR3075383} and in \cite{MR2768681}.
We are grateful to the referee of the paper for a careful reading of the paper, and for urging us
to elaborate with regard to the case of singular cardinals with countable cofinality.
We also thank Saharon Shelah for a very helpful discussion concerning several issues in pcf theory.

\newpage

\section{Pcf arrays of elementary submodels}

We define in this section the notion of a \emph{pcf array of elementary
submodels}. We provide sufficient conditions for the existence of this
kind of arrays. This will be useful in the proof of the almost strong polarized
relation later in the paper.

Let $\mu$ be a singular cardinal, $\kappa=\cf(\mu)$.
Let $\langle\mu_\alpha:\alpha<\kappa\rangle$ be an increasing sequence
of regular cardinals such that $\mu=\bigcup_{\alpha<\kappa}\mu_\alpha$.
Let $J$ be an ideal over $\kappa$ which contains the ideal $J^{\rm
  bd}_\kappa$ of all bounded subsets of $\kappa$.

The relation $f<_J g \Leftrightarrow \{\alpha\in\kappa:g(\alpha)\leq
f(\alpha)\} \in J$ defined on elements from the product
$\prod\limits_{\alpha<\kappa}\mu_\alpha$ is usually just a partial
order.
As such, it may, or may not, possess a $J$-increasing cofinal sequence of functions.
If there is a cofinal $J$-increasing sequence then the minimal length of such a
$J$-increasing sequence is denoted by ${\rm
  tcf}(\prod\limits_{\alpha<\kappa}\mu_\alpha, J)$ and called \emph{the true
cofinality} of the product.

If $J$ is a prime ideal (equivalently, if the dual of $J$ is an
ultrafilter) then the partial order defined above is actually a linear
order and hence a cofinal sequence exists. The spectrum of possible
cofinalities is the most basic concept of pcf theory.
Formally, if $\mathfrak{a}$ is a set of regular cardinals then ${\rm
  pcf}(\mathfrak{a}) = \{{\rm tcf}(\prod\mathfrak{a},\mathscr{U}):\mathscr{U}$ is an
ultrafilter over $\mathfrak{a}\}$.
For the majority of pcf theorems, $\mathfrak{a}$ has to be progressive, where a
set of regular cardinals $\mathfrak{a}$ is called \emph{progressive} iff $|\mathfrak{a}|<\min(\mathfrak{a})$.
In many cases the set $\mathfrak{a}$ is an end-segment of the set of all regular cardinals below a given
singular cardinal $\mu$.
For this set to be progressive one has to assume that $\mu$ is not a fixed point of the $\aleph$-function.

The elements of ${\rm pcf}(\mathfrak{a})$ can be characterized by ideals of the form $J_{<\lambda}[\mathfrak{a}]$.
To define these ideals, suppose that $\mathfrak{a}$ is progressive and $\mathfrak{b}\subseteq\mathfrak{a}$.
One says that $\mathfrak{b}$ \emph{dictates} cofinality less than $\lambda$ if ${\rm tcf}(\prod\mathfrak{a},\mathscr{U})<\lambda$ whenever $\mathscr{U}$ is an ultrafilter over $\mathfrak{a}$ and $\mathfrak{b}\in\mathscr{U}$.
The ideal $J_{<\lambda}[\mathfrak{a}]$ is the collection of subsets of $\mathfrak{a}$ which dictate cofinality less than $\lambda$.
By classical pcf theorems, if $\lambda\in{\rm pcf}(\mathfrak{a})$ then $J_{<\lambda^+}[\mathfrak{a}]=J_{<\lambda}[\mathfrak{a}]+\mathfrak{b}$ for some $\mathfrak{b}\subseteq\mathfrak{a}$.
Namely, the ideal $J_{<\lambda^+}[\mathfrak{a}]$ is generated over the ideal $J_{<\lambda}[\mathfrak{a}]$ by a single set.
This property is called \emph{normality}, and the set $\mathfrak{b}$ is called \emph{a generator}.
It is not unique as a set, but if $\mathfrak{b}_0,\mathfrak{b}_1$ both generate $J_{<\lambda^+}[\mathfrak{a}]$ over $J_{<\lambda}[\mathfrak{a}]$ then the symmetric difference $\mathfrak{b}_0\triangle\mathfrak{b}_1$ belongs to $J_{<\lambda}[\mathfrak{a}]$.

A fundamental theorem of pcf theory is that if $\mu>\cf(\mu)=\kappa$
then there exists a sequence
$\langle\mu_\alpha:\alpha\in\kappa\rangle$ such that ${\rm
  tcf}(\prod_{\alpha\in\kappa}\mu_\alpha,J^{\rm bd}_\kappa)=\mu^+$.
According to this theorem we can realize $\mu^+$ as a true cofinality using the
ideal of bounded subsets of $\kappa$.
In general, it might happen that $\lambda\in{\rm pcf}(\mathfrak{a})$
where $\mathfrak{a}\subseteq{\rm Reg}\cap\mu$ and $\lambda$ cannot be
expressed as a true cofinality with the ideal $J^{\rm bd}_\kappa$.
But in the case of $\lambda=\mu^{++}$ we have the following result
from \cite[Chapter VIII, Theorem 1.1]{MR1318912} at singular cardinals of uncountable cofinality.
Since the case of countable cofinality is not explicit in this monograph, we unfold the argument.
See page~\pageref{countablecof} of this paper for a more detailed explanation.
We also refer the reader to \cite{MR1356541} at this point.

\begin{theorem}
\label{thmjbd} Assume that $\mu>\cf(\mu)=\kappa\geq\aleph_0$, and
$\mu$ is a strong limit cardinal. Assume further that $2^\mu > \mu^+$.
Then there exists an increasing sequence of regular cardinals
$\langle\lambda_i: i<\kappa\rangle$ such that
$\mu=\bigcup_{i<\kappa}\lambda_i$ and ${\rm
  tcf}(\prod_{i<\kappa}\lambda_i,J^{\rm
  bd}_\kappa)=\mu^{++}$.
\end{theorem}

\par\noindent\emph{Proof}. \newline
Firstly, we argue that $\mu^{++}\leq pp(\mu)$ (there is some ideal $J$ on $\kappa$ and some increasing
sequence of regular cardinals $\langle \lambda_i:i<\kappa \rangle $ such that ${\rm
tcf}(\prod_{i<\kappa}\lambda_i,J)\ge \mu^{++}$).
Then we shall prove that this can be realized by the ideal $J^{\rm bd}_\kappa$.
For the first task, notice that $2^\mu=\mu^\kappa$ as $\mu$ is a strong limit singular cardinal.
For simplicity, we assume that $\mu$ is not a fixed point of the $\aleph$-function.\footnote{The statement of the theorem holds at fixed points of the aleph function as well, for details see \cite{MR1356541}.}
Hence there is an interval of regular cardinals $\mathfrak{a}\subseteq\mu$ such that $|\mathfrak{a}|^+<\min(\mathfrak{a})$ and even $2^{|\mathfrak{a}|}<\min(\mathfrak{a})$.
Under the above assumption, ${\rm pcf}(\mathfrak{a})$ is an interval of regular cardinals as well.
From \cite[Theorem 5.1]{MR1086455} we know that $\max{\rm pcf}(\mathfrak{a}) = |\prod\mathfrak{a}|=\mu^\kappa=2^\mu$.
Thus, $pp(\mu)\geq\max{\rm pcf}(\mathfrak{a})\geq\mu^{++}$ as we are assuming that $2^\mu>\mu^+$.

Now, since ${\rm pcf}(\mathfrak{a})$ is an interval of regular cardinals (this is a consequence of the no-holes theorem, see \cite[Theorem 3.1]{MR2768693} or \cite[Corollary 2.2]{MR1086455}) one concludes that $\mu^{++}\in{\rm pcf}(\mathfrak{a})$.
Specifically, for some increasing sequence $\langle\lambda_i:i\in\kappa\rangle$ of regular cardinals
so that $\mu=\bigcup_{i\in\kappa}\lambda_i$ it is true that $\mu^{++}={\rm
tcf}(\prod_{i\in\kappa}\lambda_i,J)$ where $J\supseteq J^{\rm bd}_\kappa$.

We move to the second task, that is, we explain why we can assume that the ideal $J$ can be taken as $J^{\rm bd}_\kappa$ itself.
For every $\lambda\in{\rm pcf}(\mathfrak{a})$ we fix a generator $\mathfrak{b}_\lambda$.
By classical pcf theorems (see, e.g., \cite[Corollary 4.4]{MR1086455}) one has $\lambda={\rm tcf}(\prod\mathfrak{b}_\lambda,J_{<\lambda}[\mathfrak{a}])$.
In particular, $\mu^{++}={\rm tcf}(\prod\mathfrak{b}_{\mu^{++}},J_{<\mu^{++}}[\mathfrak{a}])$.

Let us focus on the two generators $\mathfrak{b}_{\mu^+}$ and $\mathfrak{b}_{\mu^{++}}$.
The crucial point is that we may assume that $\mathfrak{b}_{\mu^+}\cap\mathfrak{b}_{\mu^{++}}=\varnothing$ upon replacing $\mathfrak{b}_{\mu^{++}}$ by $\mathfrak{b}_{\mu^{++}}-\mathfrak{b}_{\mu^+}$, see \cite[Theorem 4.8]{MR2768693} and \cite[Notation 4.9]{MR2768693}.
Recall that $\mu$ is a singular cardinal, and hence $J_{<\mu^+}[\mathfrak{a}]=J_{<\mu}[\mathfrak{a}]$.
But $J_{<\mu}[\mathfrak{a}]$ is none other than the ideal $J^{\rm bd}_\kappa$, since $\mu=\bigcup\mathfrak{a}$.
Finally, since $J_{<\mu^{++}}[\mathfrak{a}]=J_{<\mu^+}[\mathfrak{a}]+\mathfrak{b}_{\mu^{++}}$ and $\mathfrak{b}_{\mu^{+}}$ is disjoint from $\mathfrak{b}_{\mu^{++}}$, the ideal $J_{<\mu^{++}}[\mathfrak{a}]$ equals the ideal $J_{<\mu^+}[\mathfrak{a}]$ over the generator $\mathfrak{b}_{\mu^{++}}$.
Thus, after replacing $\mathfrak{b}_{\mu^{++}}$ by $\mathfrak{b}_{\mu^{++}}-\mathfrak{b}_{\mu^+}$ one has $J_{<\mu^{++}}[\mathfrak{a}]= J^{\rm bd}_\kappa$, and we are done.

\hfill \qedref{thmjbd}

One of the virtues of $J^{\rm bd}_\kappa$ is that if we take any
unbounded subsequence
$\langle\lambda_{i_\beta}:\beta\in\kappa\rangle$ of the above
sequence then ${\rm
  tcf}(\prod_{\beta\in\kappa}\lambda_{i_\beta},J^{\rm
  bd}_\kappa)=\mu^{++}$ as well.
Indeed, if $(f_\alpha:\alpha\in\mu^{++})$ witnesses ${\rm
  tcf}(\prod_{i<\kappa}\lambda_i,J^{\rm bd}_\kappa)=\mu^{++}$ then the restriction of the scale to the subsequence $\langle\lambda_{i_\beta}:\beta\in\kappa\rangle$ witnesses ${\rm
  tcf}(\prod_{\beta\in\kappa}\lambda_{i_\beta},J^{\rm
  bd}_\kappa)=\mu^{++}$.

In general, it is possible that $\lambda={\rm tcf}(\prod_{i<\kappa}\lambda_i,J)$ for some $J\supsetneqq J^{\rm bd}_\kappa$ and $\lambda$ cannot be realized as the true cofinality of some sequence of regular cardinals below $\mu$ with the ideal $J^{\rm bd}_\kappa$.
In such cases we do not know how to carry out our argument.
The problematic point is that by thinning-out the sequence to $\langle\lambda_{i_\beta}:\beta\in\kappa\rangle$ we may change the true cofinality.
However, in some cases one can ensure that the true cofinality is obtained by the ideal $J^{\rm bd}_\kappa$.
A notable example is when $\mu>\cf(\mu)=\kappa>\aleph_0$ and $\kappa$ is weakly compact, as proved in \cite{MR3056300}.
Another case, crucial for our proof, is the case in which $\lambda=\mu^{++}$, see \cite[Theorem 1.4]{MR1606515}.

In what follows, we will define matrices of models $(M_{\alpha
i})_{\alpha<\lambda,i<\kappa}$, where $\mu>\cf(\mu)=\kappa$ and
$\lambda\geq\mu^+$; furthermore, for some increasing sequence of regular cardinals
$\langle\lambda_i:i<\kappa\rangle$,
$\mu=\bigcup_{i<\kappa}\lambda_i$. This last part connects the
$\lambda\times\kappa$-matrix to the statement $\lambda={\rm
  tcf}(\prod_{i<\kappa}\lambda_i,J)$.


  \begin{definition}[Pcf array]
    \label{defarray}
Let $\mu$ be a singular cardinal, $\kappa=\cf(\mu)$ and $\chi$ a
sufficiently large regular cardinal above $2^\mu$. \newline
Assume that $\mu<\lambda<{\rm
tcf}(\prod_{i<\kappa}\lambda_i,J)$, $J\supseteq J_\kappa^{\rm bd}$ and for
some increasing sequence of regular cardinals $\bar{\lambda}=\langle
\lambda_i\mid
i<\kappa\rangle$, $\mu=\bigcup_{i<\kappa}\lambda_i$. \newline
A $\lambda\times\kappa$ \emph{pcf array of models} $\mathcal{M}$ is a
sequence of elementary submodels of $\mathcal{H}(\chi)$ of the form
$\langle M_{\alpha i}: \alpha<\lambda,i<\kappa\rangle$ such that (see
Figure~\ref{fig:pcfarray}):
\begin{enumerate}
\item [$(\aleph)$] For every $\alpha\in\lambda, M_{\alpha i}\prec
  M_{\alpha j}$ whenever $i<j$.
\item [$(\beth)$] For every $\alpha<\beta<\lambda$ there exists $u\in
  J$ such that $M_{\alpha j}\prec M_{\beta j}$ whenever $j\in\kappa -
  u$.
\item [$(\gimel)$] $|M_{\alpha i}|<\lambda_i$ for every
  $\alpha\in\lambda, i\in \kappa$.
\end{enumerate}
\end{definition}

We shall say that $\mathcal{M}$ is \emph{based on the sequence}
$\bar{\lambda}$ \emph{and the ideal} $J$. The following figure
demonstrates the idea in the case where $J = J^{\rm bd}_\kappa$; in
this case, for every $\alpha<\beta<\lambda$ there exists an ordinal
$i_{\alpha\beta}<\kappa$ such that $M_{\alpha j}\prec M_{\beta j}$
whenever $i_{\alpha\beta}\leq j<\kappa$.

\bigskip
\begin{center}
    \begin{figure}

\begin{tikzpicture}[font=\ttfamily,
array1/.style={matrix of nodes}]

\matrix[array1] (array1) {
	$\vdots$ & & $\vdots$ & & $\vdots$ & \reflectbox{$\ddots$}\\
$M_{\beta 0}$ & $\prec$ & $M_{\beta i_{\alpha\beta}}$ & $\prec$ & $M_{\beta i}$ & $\cdots$\\
 $\vdots$& & \rotatebox{90}{$\prec$} & & \rotatebox{90}{$\prec$} &\\
$M_{\alpha 0}$ & $\prec$ & $M_{\alpha i_{\alpha\beta}}$ & $\prec$ & $M_{\alpha
i}$ & $\cdots$\\
 $\vdots$& & $\vdots$ & & $\vdots$
 & $\cdots$\\
$M_{00}$ & $\prec$ & $M_{0 i_{\alpha\beta}}$ & $\prec$ & $M_{0 i}$ & $\cdots$\\};
\draw[->]([yshift=-3mm]array1-6-1.south west) --
([yshift=-3mm]array1-6-6.south east);

\begin{scope}[on background layer]
\fill[gray!30] (array1-2-3.north west) rectangle (array1-4-3.south east);
\fill[gray!30] (array1-2-5.north west) rectangle (array1-4-5.south east);
\end{scope}

\draw[->]([xshift=-3mm]array1-6-1.south west) --
([xshift=-5mm]array1-1-1.north west);
\node [anchor=east] at ([xshift=-6mm]array1-1-1.north west) (k) {$\lambda$};
\node [anchor=north] at ([yshift=-4mm]array1-6-6.south east) (l) {$\kappa$};
\end{tikzpicture}
\caption{A pcf-array}\label{fig:pcfarray}
\end{figure}	
\end{center}

\bigskip
A central concept related to pcf arrays is the following:

\begin{definition}[The characteristic sequence]\label{defchar}
	Let $\mathcal{M}=\langle M_{\alpha i}: \alpha<\lambda, i<\kappa\rangle$
be a pcf array of models based on $\bar{\lambda} = \langle
\lambda_i:i<\kappa \rangle$.
For every $\alpha<\lambda,i<\kappa$ let $f_\alpha(i) = \sup(M_{\alpha
  i}\cap \lambda_i)$. Each function
$f_\alpha\in\prod\limits_{i<\kappa}\lambda_i$ is called the
\emph{characteristic function} of $\mathcal{M}$ at stage $\alpha$, and
$(f_\alpha)_{\alpha<\lambda}$ is called the \emph{characteristic
  sequence}.
\end{definition}

Notice that each $f_\alpha$ is an element of
$\prod\limits_{i<\kappa}\lambda_i$ since $|M_{\alpha i}|<\lambda_i$.
The main point of the proof in the next section is that for a pcf
array $\mathcal{M}$ and $\lambda$-many functions $f_\alpha$ as above
there is a single function $h\in\prod\limits_{i<\kappa}\lambda_i$
which bounds them all.
The reason is that ${\rm
  tcf}(\prod\limits_{i<\kappa}\lambda_i,J)>\lambda$.
Now the values of $h$ will be ordinals outside the pertinent models; later in
the proof we will use the fact that any first order property of them will be
reflected to many ordinals below.

\begin{claim}
\label{clmexist} Let $\mu$ be a strong limit singular cardinal,
$\kappa=\cf(\mu)>\omega$ and $2^\mu>\mu^+$.
Let $\bar{\mu} = \langle\mu_i:i<\kappa\rangle$ be an increasing
continuous sequence of cardinals such that
$\mu=\bigcup_{i<\kappa}\mu_i$ and $\mu_0>\kappa$. Fix also another
increasing sequence of cardinals  $\bar{\lambda} =
\langle\lambda_i:i<\kappa\rangle$ such that
$\mu=\bigcup_{i<\kappa}\lambda_i$.\newline
There exists a pcf array $\mathcal{M}=\langle M_{\alpha
  i}:\alpha<\mu^+, i<\kappa\rangle$ based on $\bar{\lambda} =
\langle\lambda_i:i<\kappa\rangle$ and $J^{\rm bd}_\kappa$ such that
${}^{\mu_i^+}M_{\alpha i}\subseteq M_{\alpha i}$ for every
$\alpha\in\mu^+,i<\kappa$ and $\alpha \subseteq
\bigcup_{i<\kappa}M_{\alpha i}$ for every $\alpha\in\mu^+$.
\end{claim}

\par\noindent\emph{Proof}. \newline
We may assume, without loss of generality, that
$2^{(\mu_i^+)}<\lambda_i$ for every $i<\kappa$ by taking a subsequence
of $\bar{\lambda}$ if needed.
We use here the fact that $\mu$ is a strong limit cardinal.
Choose for every $\alpha\in\mu^+,i<\kappa$ some $M_{\alpha i}\prec
\mathcal{H}(\chi)$ such that:
\begin{enumerate}
\item [$(a)$] $|M_{\alpha i}|= 2^{(\mu_i^+)}$ and
  $2^{(\mu_i^+)}\subseteq M_{\alpha i}$.
\item [$(b)$] ${}^{\mu_i^+}M_{\alpha i}\subseteq M_{\alpha i}$.
\item [$(c)$] $\alpha\in M_{\alpha i}$ and $\alpha \subseteq
  \bigcup_{i<\kappa}M_{\alpha i}$.
\item [$(d)$] $\{M_{\alpha i}:(\alpha,i)<_{\rm lex}(\beta,j)\}\in
  M_{\beta j}$.
\end{enumerate}
The only non-trivial requirement in the construction is $\alpha
\subseteq \bigcup_{i<\kappa}M_{\alpha i}$, since typically
$\mu<\alpha$ and $|M_{\alpha i}|<\mu$ for every $i\in\kappa$. However,
$\alpha<\mu^+$ so $\alpha$ is expressible as
$\bigcup_{i<\kappa}\alpha_i$ such that $|\alpha_i|\leq 2^{(\mu_i^+)}$
for every $i\in\kappa$, and the construction follows.

We now explain why the construction works.
Requirement $(\gimel)$ of Definition \ref{defarray} is guaranteed by
$(a)$ and our assumption on the $\lambda_i$'s, so we check $(\aleph)$
and $(\beth)$.
We first note that if $(\alpha,i)<_{\rm lex}(\beta,j)$ and $\alpha\in
M_{\beta j}$ then $M_{\alpha i}\in M_{\beta j}$.
This follows from $(d)$, and from the fact that $i\in M_{\beta j}$ as
$\kappa\in 2^{(\mu_i^+)}\subseteq M_{\beta j}$.
Second, if $M_{\alpha i}\in M_{\beta j}$ and $j\leq i$ then $M_{\alpha
  i}\subseteq M_{\beta j}$.
Indeed, $\mathcal{H}(\chi)\models |M_{\alpha i}|= 2^{(\mu_i^+)}$ so
there is a function $g: 2^{(\mu_i^+)}\rightarrow M_{\alpha i}$ which
enumerates its elements. Since $M_{\beta j}\prec\mathcal{H}(\chi)$
there is such a function in $M_{\beta j}$ and since
$2^{(\mu_i^+)}\leq 2^{(\mu_j^+)}\subseteq M_{\beta j}$ we see that
the range of this function is a subset of $M_{\beta j}$, namely
$M_{\alpha i}\subseteq M_{\beta j}$.
Finally, if $M_{\alpha i}\subseteq M_{\beta j}$ then $M_{\alpha
  i}\prec M_{\beta j}$ since both are elementary submodels of the same
$\mathcal{H}(\chi)$.

Fix any $\alpha\in\mu^+$.
If $i<j$ then $(\alpha,i)<_{\rm lex}(\alpha,j)$.
Since $\alpha\in M_{\alpha j}$ by $(c)$, we infer that $M_{\alpha i}\prec
M_{\alpha j}$ so $(\aleph)$ holds.
To check $(\beth)$, assume that $\alpha<\beta$.
Since $\beta\subseteq\bigcup_{i\in\kappa}M_{\beta i}$ and
$\alpha<\beta$ there is $i_{\alpha\beta}\in\kappa$ such that $\alpha\in
M_{\beta i_{\alpha\beta}}$.
If $j\geq i_{\alpha\beta}$ then $M_{\beta i_{\alpha\beta}}\subseteq
M_{\beta j}$ and hence $\alpha\in M_{\beta j}$.
As $(\alpha,j)<_{\rm lex}(\beta,j)$ we have $M_{\alpha j}\prec
M_{\beta j}$ thus proving $(\beth)$.

\hfill \qedref{clmexist}

\newpage

\section{Almost strong relations}

In this section we shall prove the combinatorial results concerning the polarized relation.
We first provide a useful lemma which says, roughly, that if a coloring
of pairs depends only on one coordinate then it has large
monochromatic products. This observation is pertinent for both
cases of countable and uncountable cofinality.

\begin{lemma}
\label{lem0} Assume that $\mu>\cf(\mu), \tau<\mu^+$ and
$\theta<\cf(\mu)$. \newline
Let $c:\mu^+\times\mu\rightarrow\theta$ be any coloring. \newline
If $A\subseteq\mu^+, {\rm otp}(A)=\tau, B_0\in[\mu]^\mu$ and
$c\upharpoonright(A\times B_0)$ depends only on the second coordinate
then there is a set $B\in[B_0]^\mu$ such that
$c\upharpoonright(A\times B)$ is constant.
\end{lemma}

\par\noindent\emph{Proof}. \newline
Define a function $d:B_0\rightarrow\theta$ as follows.
Fix any $\alpha\in A$ and let $d(\beta)=c(\alpha,\beta)$ for every
$\beta\in B_0$. The choice of $\alpha$ is unimportant by the
assumption of the lemma.
Since $|B_0|=\mu$ and $\theta<\cf(\mu)$ one can find $j\in\theta$ and
$B\in[B_0]^\mu$ such that $\beta\in B\Rightarrow d(\beta)=j$.
It follows that $c(\alpha,\beta)=j$ for every $\alpha\in A,\beta\in
B$, so we are done.

\hfill \qedref{lem0}

The above lemma will be used within the proof of the main result of the paper, which reads as follows.

\begin{theorem}
\label{thmmt} Assume that $\mu>\cf(\mu)=\kappa\geq\aleph_0$, $\mu$ is a
strong limit cardinal and $2^\mu>\mu^+$.
Then $\binom{\mu^+}{\mu} \rightarrow \binom{\tau}{\mu}_{<\cf(\mu)}$
for every $\tau<\mu^+$.
\end{theorem}

\par\noindent\emph{Proof}. \newline
Fix two sequences of cardinals, $\bar{\mu} =
\langle\mu_i:i<\kappa\rangle$ and $\bar{\lambda} =
\langle\lambda_i:i<\kappa\rangle$, with the following properties:
\begin{enumerate}
\item [$(a)$] $\bar{\mu}$ is increasing and continuous, $\kappa<\mu_0$
  and $\mu=\bigcup_{i<\kappa}\mu_i$.
\item [$(b)$] $\bar{\lambda}$ is an increasing sequence of regular
  cardinals, $\mu=\bigcup_{i<\kappa}\lambda_i$.
\item [$(c)$] For every $i<\kappa, 2^{\mu_i^+}<\lambda_i$.
\item [$(d)$] ${\rm tcf}(\prod_{i<\kappa}\lambda_i, J^{\rm
    bd}_\kappa)=\mu^{++}$.
\end{enumerate}
\label{pageproperties}
Part $(d)$ is possible by Theorem \ref{thmjbd}.
Part $(c)$ is possible since $\mu$ is a strong limit cardinal and $\bar{\lambda}$
can be replaced by any unbounded subsequence while keeping $(d)$.
The fact that we use $J^{\rm bd}_\kappa$ in the representation of $\mu^{++}$ at $(d)$ is crucial here.

Fix any $\theta<\kappa$ and a coloring
$c:\mu^+\times\mu\rightarrow\theta$.
Let $\tau<\mu^+$ be an ordinal. By the monotonicity of the arrow
notation we may increase $\tau$ and therefore assume that
$\cf(\tau)=\kappa$ and $\tau=\bigcup_{i<\kappa}\delta_i$ where each
$\delta_i$ is indecomposable and $\delta_i>\bigcup_{j<i}\delta_j$ for
every $i<\kappa$.
We also choose another sequence of indecomposable ordinals $\langle\nu_i:i\in\kappa\rangle$ so that $\delta_i<\nu_i<\delta_{i+1}$ for every $i\in\kappa$.
These assumptions are possible since $S^{\mu^+}_\kappa$ is stationary
and the set of indecomposable ordinals below $\mu^+$ is a club.

Let $\chi$ be a sufficiently large regular cardinal in which the arguments of our proof will be carries out.
Let $\mathcal{M} = \langle M_{\alpha i}:\alpha<\mu^+, i<\kappa\rangle$
be a pcf array of submodels of $\mathcal{H}(\chi)$, based on the ideal
$J^{\rm bd}_\kappa$ and the sequences $\bar{\mu},\bar{\lambda}$ and let
$\langle f_\alpha: \alpha<\mu^+\rangle $ be its characteristic sequence.
We may assume that ${}^\kappa\theta\subseteq M_{\alpha i}$ for every $\alpha\in\mu^+,i\in\kappa$.
Since ${\rm tcf}(\prod_{i<\kappa}\lambda_i, J^{\rm
bd}_\kappa)=\mu^{++}$, we can fix a function\label{fstarpage}
$f^*\in\prod_{i<\kappa}\lambda_i$ such that $\alpha\in\mu^+
\Rightarrow f_\alpha<_{J^{\rm bd}_\kappa}f^*$.

We shall use $f^*$ as a \emph{translation} of the coloring $c$, defined
on pairs, to a function defined on singletons in the following sense: for every
$\alpha\in\mu^+$ define
$g_\alpha\in{}^\kappa\theta$ by letting $g_\alpha(i) =
c(\alpha,f^*(i))$. Here we still depend on both coordinates, but since
$\theta^\kappa<\mu$ there is a fixed $g^*\in{}^\kappa\theta$ such that
the set $H = \{\alpha\in\mu^+:g_\alpha=g^*\}$ is of size $\mu^+$.
Now $c(\alpha,f^*(i))$ depends only on the righthand coordinate when
we restrict ourselves to $H$ in the lefthand coordinate: if
$\alpha_1,\alpha_2\in H$ then $c(\alpha_1,f^*(i))=g_{\alpha_1}(i)=g^*(i)
=g_{\alpha_2}(i)=c(\alpha_2,f^*(i))$.

The role of $\mathcal{M}$ and $f^*$ is to create the small component
of the monochromatic product. For the large component we choose a
continuous increasing sequence $\langle N_\zeta:\zeta\in\mu^+\rangle$
of submodels of $\mathcal{H}(\chi)$ so that for every $\zeta\in\mu^+$
the following requirements are met:
\begin{enumerate}
\item [$(\aleph)$] $\mu\subseteq N_\zeta, |N_\zeta|=\mu$.
\item [$(\beth)$] $c, H, g^*, \mathcal{M} \in N_\zeta$.
\item [$(\gimel)$] $N_\zeta\in N_{\zeta+1}$.
\end{enumerate}
Notice that the set $E = \{\zeta\in\mu^+:N_\zeta\cap\mu^+=\zeta\}$ is
a club subset of $\mu^+$.
It would be helpful to concentrate on members of $E$.
Since $E$ is a club we may assume, without loss of generality, that
every $\delta_i$ belongs to ${\rm acc}(E)$.
For every $i<\kappa$ we choose an ordinal $\alpha(i)\in\delta_i$ such
that $\alpha(i)>\bigcup\{\nu_j:j<i\}$.
Also, pick $\alpha(*)\in H - \tau$.

The monochromatic product will be created now by induction.
For every $i<\kappa$ we build sets $A_i,B_i$ which will approximate
the large and the small components of the product respectively.
Likewise, we choose an ordinal $j(i)\in\kappa$ and we try to keep the
following requirements:\label{reqinductpage}
\begin{enumerate}
\item [$(a)$] $j(i)>i$; also, $i_0<i_1\Rightarrow
  \lambda_{j(i_0)}<\mu_{j(i_1)}$.
\item [$(b)$] If $\sigma_0,\sigma_1\in \{\delta_\ell:\ell\leq i\} \cup
  \{\alpha(\ell):\ell\leq i\} \cup \{\alpha(*)\}$ and
  $\sigma_0<\sigma_1$ then $f_{\sigma_0}(j)<f_{\sigma_1}(j)$ for every
  $j\in[j(i),\kappa)$.
\item [$(c)$] $A_i\subseteq H\cap\delta_i, {\rm otp}(A_i)=\delta_i$
  and $A_i\in M_{\nu_i j(i)}$.
\item [$(d)$] $B_i\subseteq\lambda_{j(i)} -
  \bigcup\{\lambda_{j(\ell)}:\ell<i\}, {\rm otp}(B_i)=\lambda_{j(i)}$
  and $B_\ell\in M_{\nu_i j(i)}$ for every $\ell<i$.
\item [$(e)$] If $\alpha\in\bigcup\{A_\ell:\ell\leq
  i\}\cup\{\alpha(*)\}$ and $\beta\in B_\ell$ for some $\ell\leq i$
  then $c(\alpha,\beta) = g^*(j(\ell))$.
\end{enumerate}
Suppose that we can carry the induction.
Define $A = \bigcup_{i<\kappa}A_i\cup\{\alpha(*)\}, B = \bigcup_{i<\kappa}B_i$.
Notice that $|B| = \bigcup\{\lambda_{j(i)}:i<\kappa\}=\mu$ and ${\rm otp}(A) =
\tau+1$.
Choose any pair of ordinals $\alpha\in A, \beta\in B$.
Let $i<\kappa$ be the first ordinal so that $\beta<\lambda_{j(i)}$ and
$\alpha\in\bigcup\{A_\ell:\ell\leq i\}\cup\{\alpha(*)\}$.
Requirement $(d)$ implies that the $B_\ell$s are mutually disjoint.
By the disjointness of the $B_\ell$s there is a unique $\ell\leq i$ such that
$\beta\in B_\ell$, and then $c(\alpha,\beta)=g^*(j(\ell))$.
This means that $c(\alpha,\beta)$ depends only on the righthand coordinate (as
$A\subseteq H$), so by Lemma \ref{lem0} we are done.

It remains to carry out the induction.
Suppose that $i<\kappa$, and assume that $A_\ell,B_\ell,j(\ell)$ have been
defined for every $\ell<i$.
We shall choose $j(i)$ and build $A_i$ together, and then we shall describe the
construction of $B_i$.
For every $\ell<i$ we pick an ordinal $\eta_\ell\in\kappa$ such that
$A_\ell,B_\ell,j(\ell)\in M_{\alpha(i)\eta_\ell}$. Notice that for every
$j\in[\eta_\ell,\kappa)$ we will have $A_\ell,B_\ell,j(\ell)\in M_{\alpha(i)j}$
  by the properties of the array $\mathcal{M}$.
The choice of $\eta_\ell$ is possible since $\nu_\ell<\alpha(i)$ and hence
there exists $\eta_\ell\in\kappa$ such that $M_{\nu_\ell j}\prec
M_{\alpha(i)j}$ for every $j\in[\eta_\ell,\kappa)$.
Now $A_\ell, B_\ell, j(\ell)\in M_{\nu_\ell j(\ell)}$ by the induction
hypothesis, so they belong to $M_{\nu_\ell j}$ for every $j\geq j(\ell)$.
It follows that for every sufficiently large $j<\kappa$ we have $A_\ell,
B_\ell, j(\ell)\in M_{\alpha(i)j}$.
The above argument works separately for every $\ell<i$, and we can choose now
$j_0<\kappa$ such that $A_\ell, B_\ell, j(\ell)\in M_{\alpha(i) j_0}$
simultaneously for every $\ell<i$.
This is possible since $i<\kappa$.

We choose an increasing sequence $\langle
\zeta_\varepsilon:\varepsilon<\cf(\delta_i)\rangle$ of elements of $E$ such
that $\zeta_0>\alpha(i)$ and $\delta_i = \bigcup\{\zeta_\varepsilon:
\varepsilon<\cf(\delta_i)\}$.
Notice that $M_{\alpha(i)j_0}\prec N_{\zeta_\varepsilon}$ for every
$\varepsilon<\cf(\delta_i)$:
fixing $\varepsilon<\cf(\delta_i)$, we
first note that $M_{\alpha(i)j_0}\in N_{\zeta_\varepsilon}$,
using that both ordinals $\alpha(i)$ and $j_0$ belong to
$N_{\zeta_\varepsilon}$, as $\mathcal{M}\in N_{\zeta_\varepsilon}$ and
$j_0\in \kappa\subseteq N_{\zeta_\varepsilon}$. Now, since
$|M_{\alpha(i)j_0}|<\lambda_{j_0}$ (by the definition of a pcf array), we have
in
particular $|M_{\alpha(i)j_0}|<\mu$; since $\mu\subseteq
N_{\zeta_\varepsilon}$, $M_{\alpha(i)j_0}\subseteq
N_{\zeta_\varepsilon}$. Finally, as both  $M_{\alpha(i)j_0}$ and
$N_{\zeta_\varepsilon}$ are elementary submodels of ${\mathcal
H}(\chi)$ we may conclude that $M_{\alpha(i)j_0}\prec
N_{\zeta_\varepsilon}$.

Fix any $\varepsilon<\cf(\delta_i)$.
We define a first order formula $\varphi(x)$ with parameters ($x$ is the free
variable), as follows:
$$
x\in H \wedge x>\zeta_\varepsilon \wedge (\forall\ell<i) (\beta\in B_\ell
\rightarrow c(x,\beta)=g^*(j(\ell))).
$$
Remark that $\mathcal{H}(\chi)\models\varphi (\alpha(*))$ using the induction
hypothesis. Notice that $\varphi(x)$ is definable in
$N_{\zeta_{\varepsilon+1}}$, so by elementarity there is a set
$\Gamma_\varepsilon\subseteq N_{\zeta_{\varepsilon+1}} -
N_{\zeta_{\varepsilon}}$ such that ${\rm otp}(\Gamma_\varepsilon) =
\zeta_{\varepsilon+1}$ and $\gamma\in \Gamma_\varepsilon \Rightarrow
\mathcal{H}(\chi)\models \varphi(\gamma)$.
Let $A_i = \bigcup\{\Gamma_\varepsilon:\varepsilon<\cf(\delta_i)\}$.
Observe that ${\rm otp}(A_i)=\delta_i$.

We can choose now $j(i)\in\kappa$ by setting the following requirements.
First we find $j_1\in\kappa$ such that $A_i\subseteq M_{\nu_ij_1}$ and then we choose $j(i)\geq\max\{j_0,j_1,i\}$.
Now we increase $j(i)$ if needed, by considering a few more requirements.
We make sure that $A_i\in M_{\nu_ij_1}$.
This can be done since $\delta_i<\nu_i$ and hence every subset of $\delta_i$ is bounded in $\bigcup_{j\in\kappa}M_{\nu_i j}$.
In particular, $A_i\subseteq\delta_i$, so for some $j_1\in\kappa$ we see that $A_i\subseteq M_{\nu_ij_1} \subseteq M_{\nu_ij(i)}$ and by increasing the above chosen $j(i)$ we arrive at a sufficiently closed submodel so that inclusion implies membership.
If still needed we increase $j(i)$ once more to get $M_{\nu_ij(i)} \prec M_{\alpha(*)j(i)}$ and $f_{\delta_i}(j)<f^*(j)$ for every $j\geq j(i)$.

Finally, we define the set $B_i$.
Let $\Gamma = \{\beta\in\lambda_{j(i)}: \forall\alpha\in \bigcup_{\ell\leq i}A_\ell \cup\{\alpha(*)\}, c(\alpha,\beta)=g^*(j(i))\}$.
Notice that $\Gamma$ is definable in $M_{\alpha(*)j(i)}$, since all the parameters are in this model.
This is true, in particular, for $g^*$, since all the functions from $\kappa$ to $\theta$ are in every $M_{\alpha i}$.
The focal point here is that $f^*(j(i))\in\Gamma$.
Indeed, $f^*(j(i))<\lambda_{j(i)}$ and for every $\alpha\in H$ we have $c(\alpha,f^*(j(i)))=g^*(j(i))$.
This is true, in particular, for every $\alpha\in \bigcup_{\ell\leq i}A_\ell \cup\{\alpha(*)\}$ as all these ordinals are from $H$.

Since $f^*(j(i))>f_{\nu_i}(j(i)) = \sup(M_{\nu_ij(i)}\cap\lambda_{j(i)})$ we infer that $f^*(j(i))\notin M_{\nu_ij(i)}$.
By elementarity there is an unbounded set of ordinals below $\lambda_{j(i)}$ in $M_{\nu_ij(i)}$ which belong to $\Gamma$.
Truncate the part of this set below $\bigcup\{\lambda_{j(\ell)}: \ell<i\}$, and let $B_i$ be the remainder.
One can verify that $A_i,B_i$ and $j(i)$ satisfy the requirements of the inductive process, so we are done.

\hfill \qedref{thmmt}

Let us add a clarification concerning singular cardinals of countable
cofinality.\label{countablecof}
  The representation theorem 1.1 is phrased in Shelah's book [12] with the assumption that $\mu>\cf(\mu)>\aleph_0$.
  In order to incorporate singular cardinals of countable cofinality one can argue in two different ways.
  Let us describe both strategies.

  Firstly, one can try to improve the representation theorem in such a way that singular cardinals with countable cofinality are included.
  Secondly, one can argue that in the specific case of countable cofinality there is no need to work with $J^{\rm bd}_\omega$, and every ideal $J\supseteq J^{\rm bd}_\omega$ will give the desired result.
  Let us examine both alternatives.

  For the first alternative, when the monograph of Shelah was published the proof of Theorem 1.1 required uncountable cofinality.
  Later, Shelah improved the methods of proof and obtained the same result at singular cardinals of countable cofinality as well.
  Thus, in [13, Theorem 1.4] Shelah indicated explicitly that a representation of $\mu^{++}$ with $J^{\rm bd}_\omega$ applies to $\mu>\cf(\mu)=\omega$ as well.
  We mention the fact that even the classical representation of $\mu^+$ at singular cardinals with countable cofinality is more involved, see [1, Exercise 2.25].
  However, for $\mu^{++}$ (and even a bit more) one can still prove the same statement for all possible cofinalities, including singular cardinals with countable cofinality.

  The second alternative is not used in the current paper, but it is worth mentioning this direction since it might be helpful in other statements of this type.
  It is based on the fact that we need the specific ideal $J^{\rm bd}_\kappa$ only at one point of the proof, where $\kappa=\cf(\mu)$.
  We deal with a sequence $\langle\lambda_i:i<\kappa\rangle$ so that $\mu^{++}={\rm tcf}(\prod_{i<\kappa}\lambda_i,J^{\rm bd}_\kappa)$, but then we may thin-out the sequence and we must be sure that the true cofinality remains $\mu^{++}$.
  An arbitrary ideal does not necessarily preserve this property, while $J^{\rm bd}_\kappa$ has this virtue.
  However, one needs only the fact that $J^{\rm bd}_\kappa$ is $\kappa$-complete, and every $\kappa$-complete ideal over $\kappa$ which extends $J^{\rm bd}_\kappa$ will have the same effect.

  In the specific case of $\kappa=\aleph_0$ we can modify the proof using any ideal $J\subseteq J^{\rm bd}_\omega$, since every such ideal is $\aleph_0$-complete (this is a special feature of $\aleph_0$, of course).
  The changes in the proof are as follows.
  Rather than choosing $j(i)\in\kappa$ one has to choose a set $u_i\in{J}$ for every $i<\kappa$ and then let $j(i)=\min(\kappa-u_i)$.
  In all the arguments which apply to an end-segment (i.e., for every $j\in[j(i),\kappa)$) one should focus on every $j\in\kappa-u_i$.
  In the construction of the elementary submodels one should require $u_i\in M_{\alpha(i)\eta_\ell}$ for each $i<\kappa$.
  Finally, instead of choosing a sufficiently large $j$ by removing less than $\kappa$ many bounded sets, one should pick an ordinal $j\in\kappa-\bigcup_{\ell\in{i}}u_\ell$.
  This is possible by the $\kappa$-completeness of $J$ and the fact that $u_\ell\in{J}$ for each $\ell\in{i}$.

  In the original form of the paper we argued in this way, and the proof was somewhat cumbersome.
  After a helpful discussion with Saharon, we realized that there is no need to argue in that way, so we simplified the proof and indicated that the required representation theorem works at singular cardinals of countable cofinality as well.
  However, it might be useful in the future to bear in mind that the argument can be based on $\kappa$-complete ideals in general, and it is not necessary to employ $J^{\rm bd}_\kappa$.

In this paper we focused on the so-called \emph{balanced} polarized relation, in which the required size of the monochromatic product is identical in both colors.
Let us add a few words with regard to the unbalanced relation.
The central assumption in our paper, as well as in \cite{MR1606515}, is $2^\mu>\mu^+$.
However, if $\binom{\mu^+}{\mu}\rightarrow\binom{\tau}{\mu}_\theta$ for every $\tau\in\mu^+$ and one collapses $2^\mu$ to $\mu^+$ then this instance of an almost strong relation is preserved by virtue of the completeness of the collapse.

Consider now the parallel unbalanced relation $\binom{\mu^+}{\mu}\rightarrow\binom{\mu^+\quad\tau}{\mu\quad\mu}$.
It was shown in \cite{MR4245063} that if $\mu>\cf(\mu)=\omega$ and $2^\mu=\mu^+$ then $\binom{\mu^+}{\mu}\nrightarrow\binom{\mu^+\quad\omega_2}{\mu\quad\mu}$, and it was shown in \cite{MR4211881} that consistently $\binom{\mu^+}{\mu}\nrightarrow\binom{\mu^+\quad\omega_1}{\mu\quad\mu}$.
For singular cardinals of countable cofinality this negative relation is optimal, since for such cardinals $\binom{\mu^+}{\mu}\rightarrow\binom{\mu^+\quad\tau}{\mu\quad\mu}$ for every $\tau\in\omega_1$, as proved in \cite{MR2367118}.
If $\cf(\mu)>\omega$ then the negative relation is even sharper, as $\binom{\mu^+}{\mu}\nrightarrow\binom{\mu^+\quad\omega}{\mu\quad\mu}$, a result of Erd\H{o}s, Hajnal and Rado.
Thus, assuming $2^\mu=\mu^+$ there is a meaningful difference between the balanced and the unbalanced almost strong polarized relations.
This is reflected by the fact that at the balanced relation one can obtain $\binom{\mu^+}{\mu}\rightarrow\binom{\tau}{\mu}_\theta$ for every $\tau\in\mu^+, \theta\in\cf(\mu)$, and this is optimal.
Indeed, $2^\mu=\mu^+$ implies the failure of the strong relation, that is $\binom{\mu^+}{\mu}\nrightarrow\binom{\mu^+}{\mu}$.

We conclude the paper with the following natural problem:

\begin{question}
  \label{qstrongrelation} Let $\mu$ be a strong limit singular cardinal and suppose that $2^\mu>\mu^+$.
  Is it consistent that $\binom{\mu^+}{\mu}\nrightarrow\binom{\mu^+}{\mu}$?
\end{question}

\bibliographystyle{alpha}

\end{document}